# THE SIZE OF COMPONENTS IN CONTINUUM NEAREST-NEIGHBOR GRAPHS

By Iva Kozakova, Ronald Meester and Seema Nanda

*Vrije Universiteit Amsterdam, Vrije Universiteit Amsterdam and University of Tennessee*

We study the size of connected components of random nearest-neighbor graphs with vertex set the points of a homogeneous Poisson point process in $\mathbb{R}^d$. The connectivity function is shown to decay superexponentially, and we identify the exact exponent. From this we also obtain the decay rate of the maximal number of points of a path through the origin. We define the *generation number* of a point in a component and establish its asymptotic distribution as the dimension $d$ tends to infinity.

**1. Basic definitions and results.** Let $X$ be a homogeneous density 1 Poisson process in $\mathbb{R}^d$ with an "extra" Poisson point located at the origin. Let $\vec{\mathcal{G}}_d$ denote the directed graph whose vertices are the Poisson points and in which there is a directed edge from $s \in X$ to $s' \in X$ if $s'$ is the nearest neighbor (NN) of $s$. Ignoring the directions of the edges leads to an undirected graph which we denote by just $\mathcal{G}_d$.

The graph $\mathcal{G}_d$ was introduced and studied in [1]. They showed that $\mathcal{G}_d$ contains a.s. no unbounded component, in any dimension $d$. In the current paper, we are interested in the tail behavior of components in $\mathcal{G}_d$; rather than only stating that they are finite, we would like to know how large these clusters typically are. To this end, we make a few definitions.

We denote by $p_d(n, L)$ the probability that there is a directed path in $\vec{\mathcal{G}}_d$ starting at the origin, touching exactly $n$ distinct points (besides the origin) and ending at a point $s$ with $|s| > L$ (where $|\cdot|$ denotes Euclidean distance). Furthermore, $\tau_d(L)$ is the probability that there is a path in $\mathcal{G}_d$ starting at the origin and ending at a point $s$ with $|s| > L$.

Since in a Poisson process the distances between any two pairs of points are different a.s., each point has an a.s. unique NN. It is quite possible that two Poisson points are each other's NN. In this case there is a mini-loop of two directed edges in $\vec{\mathcal{G}}_d$ between these two NN's. In fact, since each Poisson point must have an NN, and all components are finite, each component contains at least one such closed mini-loop. This would happen when a Poisson point $s$ is an NN of one or more other Poisson points, and when one of the latter is an NN of $s$ as well. Furthermore, each component contains *exactly* one closed mini-loop, as the existence of more than one such closed mini-loop would imply the existence of more than one NN for some Poisson point in the component.

It is not possible to have any type of closed circuit other than the one described here, as the existence of any closed circuit involving more than two Poisson points would contradict the fact that the lengths of successive directed edges are decreasing. The existence of an (a.s.) unique NN implies that along any path in $\mathcal{G}_d$ there is at most one change of direction (in the arrows of $\vec{\mathcal{G}}_d$), and the change of direction must take place at one of the two Poisson points in a loop.

We finally mention that it is a standard fact that there is a uniform and finite upper bound, depending on the dimension $d$, for the number of points that have the same point as their nearest neighbor. The maximum number of such points is called the *kissing number* and denoted by $K_d < \infty$; see, for instance, [3].

There are at least two ways to measure the size of a component in $\mathcal{G}_d$: one can look at the *diameter* of the component, or at the *number of points* in the component. First we state a result that tells us that the diameter of a component decays faster than exponentially, and which specifies the exponent exactly. For a related result in a discrete setting, see [2].

THEOREM 1.1. *There exist constants $C_1$, $C_2$, $L_0 \in (0, \infty)$ (depending only on the dimension $d$) such that*

$$(1) \qquad e^{-C_1 L (\log L)^{(d-1)/d}} \leq \tau_d(L) \leq e^{-C_2 L (\log L)^{(d-1)/d}} \qquad \text{for } L \geq L_0.$$

Turning to the number of points in a component, we let $\rho_d(n)$ denote the probability that there is a path in $\vec{\mathcal{G}}_d$ through the origin touching more than $n$ distinct points.

THEOREM 1.2. *There exist constants $C'_1$, $C'_2$, $N_0 \in (0, \infty)$ (depending only on $d$) such that*

$$e^{-C'_1 n \log n} \leq \rho_d(n) \leq e^{-C'_2 n \log n} \qquad \text{for } n \geq N_0.$$

For fixed dimension, it seems difficult to make more precise statements about the number of points in a component. Therefore, we investigate what happens in the limit when the dimension $d \to \infty$, and we will obtain some indirect information about the size of a component via the so-called *generation number* of a point. We already remarked that any component contains exactly one mini-loop of two arrows. The points in this mini-loop are given generation number 1. A Poisson point $x$ receives generation number $k$ if the graph distance to the unique mini-loop in the component of $x$ is equal to $k-1$. We denote by $g_d(k)$ the probability that the origin has generation number $k$, in dimension $d$. The following result is an indication that components typically are very small.

THEOREM 1.3. *We have that*
$$\lim_{d \to \infty} g_d(k) = \frac{k}{(k+1)!}, \qquad k = 1, 2, \ldots.$$

Section 2 contains the proofs of Theorems 1.1 and 1.2, while Section 3 contains the proof of Theorem 1.3.

**2. Proofs of Theorems 1.1 and 1.2.** In this section $W_1, \ldots, W_n$ are i.i.d. $\mathbb{R}^d$-valued random variables whose common probability density is given by $e^{-V(B(0,|w|))} = e^{-\pi_d |w|^d}$. The lower bound is based on the following. Note that $\tau_d(L) \geq p_d(n, L)$ for every $n$.

PROPOSITION 2.1. *There exist constants $b_1, c_1 \in (0, \infty)$ (depending only on $d$) such that*

$$(2) \qquad p_d(n, L) \geq \frac{(b_1)^n}{n!} e^{-c_1 n (L/n)^d}.$$

As will be seen from the proof of the proposition, one may choose any $\theta \in (0, \frac{\pi}{4}]$ and take $b_1$ to be the probability that a uniformly distributed random point on the unit sphere $S^{d-1}$ falls into the "polar cap" of opening half-angle $\theta$. The corresponding $c_1$ may then be taken as $\frac{\pi_d}{(\cos \theta)^d}$ where $\pi_d = V(B(0,1))$, where $V$ denotes Euclidean volume and $B(x,r)$ denotes the (open) ball of radius $r$ centered at $x \in \mathbb{R}^d$.

For the upper bound we need to bound not only $p_d(n, L)$ but also some closely related quantities that we now define.

For $j \in \{0, \ldots, n\}$ we define $p_d(n, L, j)$ to be the probability of the event $E(n, L, j)$, that there are two directed paths in $\vec{\mathcal{G}}_d$: one from $0$ to some $s'$, touching exactly $j$ particles (besides $0$) and one from some $s$ to the same $s'$, touching exactly $n - j$ particles (besides $s'$) and such that $|s| > L$. Thus $p_d(n, L) = p_d(n, L, n)$ and furthermore (by the properties of the directed graph $\vec{\mathcal{G}}_d$ mentioned above and) since $\tau_d(L)$ is equal to the probability of $\bigcup_{n=1}^{\infty} \bigcup_{j=1}^{n} E(n, L, j)$, we have

$$\tau_d(L) \leq \sum_{n=1}^{\infty} \sum_{j=0}^{n} p_d(n, L, j). \tag{3}$$

The upper bound is based on the following:

PROPOSITION 2.2. *There exists a constant $c_2 \in (0, \infty)$ (depending only on d) such that*

$$\begin{aligned} p_d(n, L, j) &\leq \frac{(K_d)^n}{j!(n-j)!} P(|W_1 + \cdots + W_n| \geq c_2 L) \\ &\leq \frac{(K_d)^n}{j!(n-j)!} P(|W_1| + \cdots + |W_n| \geq c_2 L), \end{aligned} \tag{4}$$

*where $W_1, \ldots, W_n$ are i.i.d. $\mathbb{R}^d$-valued random variables as described earlier.*

As a corollary to Proposition 2.2 and inequality (3), we have the following.

PROPOSITION 2.3. $\tau_d(L) \leq e^{2K} P(\mathcal{U} \geq c_2 L)$, *where $\mathcal{U}$ is a random variable with $E(e^{r\mathcal{U}}) = e^{2KE(e^{r|W_1|}-1)}$.*

PROOF OF PROPOSITION 2.1. Given $(x_1, \ldots, x_n) \in (\mathbb{R}^d)^n$ we define $s_0 = 0$ and $s_i = x_1 + \cdots + x_i$ for $i = 1, \ldots, n$. Let the set of points $(x_1, \ldots, x_n) \in (\mathbb{R}^d)^n$ satisfying the following three conditions be denoted by $\mathcal{S}$:

(i) $|x_1| \geq |x_2| \geq \cdots \geq |x_n|$,
(ii) $|s_n| \geq L$,
(iii) $s_i \notin \bigcup_{j=1}^{i} B_j$ for $i = 1, \ldots, n$, where $B_j = B(s_{j-1}, |x_j|)$ is the open ball centered at $s_{j-1}$ of radius $|x_j|$.

Note that [because of (i)] condition (iii) may be replaced by:

(iii′) $B_l \cap \{s_0, \ldots, s_n\} = s_{l-1}$ for $l = 1, \ldots, n$.

We now claim that

$$p_d(n, L) = \int_{\mathcal{S}} e^{-V(\bigcup_{j=1}^{n} B_j)} \, dx_1 \cdots dx_n. \tag{5}$$

Since we use this type of equality (a variation on Campbell's theorem) a number of times, and since we have not been able to find a proof in the literature, we spend a few lines on the proof of (5). Our proof proceeds by a suitable discretization of $\mathbb{R}^d$. For $k = 1, 2, \ldots$, we consider (nested) subdivisions of $\mathbb{R}^d$ into $d$-dimensional cubes of side length $2^{-k}$. We denote by $\mathcal{D}(\mathcal{S})$ the collection of $n$-tuples $(D_1^k, \ldots, D_n^k)$ of these cubes with the property that for all $s_1 \in D_1^k, \ldots, s_n \in D_n^k$ we have $(s_1, \ldots, s_n) \in \mathcal{S}$.

Denoting the event in question by $E$, we denote by $E_k$ the event that in addition, all points $s_1, \ldots, s_n$ of the directed path are the only Poisson points in their respective cubes $D_1^k, \ldots, D_n^k$, of side length $2^{-k}$ and such that $(D_1^k, \ldots, D_n^k) \in \mathcal{D}(\mathcal{S})$. It is clear that $E_k \to E$ and that $P(E_k) \to P(E)$, as $k \to \infty$. We therefore need to compute $\lim_{k \to \infty} P(E_k)$. Denoting the integrand in (5) by $f$, we write $\bar{f}(D_1^k, \ldots, D_n^k)$ for the expectation of $f(X_1, \ldots, X_n)$ where the $X_i$'s are independent and uniform over $D_i$, respectively, $i = 1, \ldots, n$.

Note that $E_k$ occurs when there is an $n$-tuple $(D_1^k, \ldots, D_n^k) \in \mathcal{D}(\mathcal{S})$ such that each box $D_i^k$ of side length $2^{-k}$ in this sequence contains exactly one Poisson point, such that in addition, the correct balls around these points contain no further Poisson point. This gives

$$P(E_k) = \sum_{(D_1^k, \ldots, D_n^k) \in \mathcal{D}(\mathcal{S})} (2^{-dk} + o(2^{-dk}))^n \bar{f}(D_1^k, \ldots, D_n^k)$$

$$= 2^{-dnk} \sum_{(D_1^k, \ldots, D_n^k) \in \mathcal{D}(\mathcal{S})} \bar{f}(D_1^k, \ldots, D_n^k) + o(2^{-dnk})$$

$$\to \int_{\mathcal{S}} f(x_1, \ldots, x_n) \, dx_1 \cdots dx_n,$$

as $k \to \infty$, proving (5).

Now let $\mathcal{C}(\theta)$ be the polar cap of half angle $\theta \leq \frac{\pi}{4}$ in the unit sphere of $\mathbb{R}^d$ with vertex at the origin; that is,

$$\mathcal{C}(\theta) = \{x = (x_1, \ldots, x_d) \in \mathbb{R}^d : |x| = 1 \text{ and } x_1 \geq \cos \theta\}.$$

Since

(6) $$V\left(\bigcup_{j=1}^n B_j\right) \leq \sum_{j=1}^n V(B_j),$$

it follows [using the notation of Proposition 2.2, letting $S_i = W_1 + \cdots + W_i$ and $B_j = B(S_{j-1}, |W_j|)$] that

$$p_d(n, L) \geq \int_{\mathcal{S}} \prod_{j=1}^n e^{-\pi_d |x_j|^d} \, dx_1 \cdots dx_n$$

$$= P\left(|W_1| \geq \cdots \geq |W_n|, \left|\sum_{j=1}^n W_j\right| \geq L \text{ and } S_i \notin \bigcup_{j=1}^i B_j \text{ for each } i\right)$$

$$\geq P\left(|W_1| \geq \cdots \geq |W_n|, \left|\sum_{j=1}^n W_j\right| \geq L \text{ and } \frac{W_i}{|W_i|} \in \mathcal{C}(\theta) \text{ for each } i\right)$$

$$\geq (b(\theta))^n P\left(|W_1| \geq \cdots \geq |W_n|, \sum_{j=1}^n |W_j|\cos(\theta) \geq L\right)$$

$$= \frac{(b(\theta))^n}{n!} P\left(\sum_{j=1}^n |W_j| \geq \frac{L}{\cos\theta}\right) \geq \frac{(b(\theta))^n}{n!}\left[P\left(|W_1| \geq \frac{L}{n\cos(\theta)}\right)\right]^n$$

$$= \frac{(b(\theta))^n}{n!} e^{-(\pi_d/(\cos\theta)^d)n(L/n)^d}.$$

Here $b(\theta)$ represents the probability that $S_i$ lies in a cone of half angle $\theta \leq \frac{\pi}{4}$ with vertex at $S_{i-1}$. The result now follows from the approximation $\frac{b_d(\theta)^n}{n!} \approx e^{n(\log(b(\theta))-\log n+1)}$. □

To prove Proposition 2.2 we need Lemma 2.4 below which should be compared to (6). As a replacement for (5), one can straightforwardly obtain the inequality

$$(7) \qquad p_d(n,L,j) \leq \int_{\mathcal{S}_j} e^{-V(\bigcup_{i=1}^j B_i \cup \bigcup_{i=j+1}^n B'_i)} \, dx_1 \cdots dx_n.$$

Here $B'_i = B(s_i, |x_i|)$ and $\mathcal{S}_j$ is the set of points $(x_1,\ldots,x_n) \in (\mathbb{R}^d)^n$ satisfying the three conditions (where again $s_0 = 0$ and $s_i = x_1 + \cdots + x_i$):

(i) $|x_1| \geq \cdots \geq |x_j|$ and $|x_{j+1}| \leq \cdots \leq |x_n|$,
(ii) $|s_n| \geq L$,
(iii) $B_l \cap \{s_0,\ldots,s_n\} = s_{l-1}$ for $l = 1,\ldots,j$ and $B'_l \cap \{s_0,\ldots,s_n\} = s_l$ for $l = j+1,\ldots,n$.

LEMMA 2.4. *Let* $(x_1,\ldots,x_n) \in \mathcal{S}_j$. *Then*

$$(8) \qquad V\left(\bigcup_{i=1}^j B_i \cup \bigcup_{i=j+1}^n B'_i\right) \geq \frac{1}{K_d}\left[\sum_{i=1}^j V(B_i) + \sum_{i=j+1}^n V(B'_i)\right].$$

PROOF. This is a consequence of the fact that no point $x$ can lie in the intersection of more than $K_d$ of the balls $\tilde{B}_i$. To see this, note that if this were not so, then $x$ would be the closest point (among $\{x, q_1, \ldots, q_n\}$) to more than $K_d$ of the $q_i$'s, contradicting the definition of the kissing number $K_d$. □

PROOF OF PROPOSITION 2.2. From now on we call $K_d$ merely $K$. Starting with (7), we apply Lemma 2.4 to find

$$p_d(n,L,j) \leq \int_{\mathcal{S}_j} e^{-(1/K)\sum_{i=1}^j V(B_i)} e^{-(1/K)\sum_{i=j+1}^n V(B'_i)} \, dx_1 \cdots dx_n$$

$$= \int_{\mathcal{S}_j} e^{-(1/K)\pi_d \sum_{i=1}^{j} |x_i|^d} e^{-(1/K)\pi_d \sum_{i=j+1}^{n} |x_i|^d} \, dx_1 \cdots dx_n$$

$$= \int_{\mathcal{S}_j} \prod_{j=1}^{n} e^{-(1/K)\pi_d |x_j|^d} \, dx_1 \cdots dx_n.$$

By the change of variable $x_i = y_i K^{(1/d)}$ for each $i$, the last integral becomes

(9) $$K^n \int_{\mathcal{S}'_j} \prod_{j=1}^{n} e^{-\pi_d |y_j|^d} \, dy_1 \cdots dy_n,$$

where $\mathcal{S}'_j$ is the set of $(y_1, \ldots, y_n) \in (\mathbb{R}^d)^n$ such that $(K^{(1/d)}y_1, \ldots, K^{(1/d)}y_n) \in \mathcal{S}_j$. To get an upper bound on (9) we simply drop the third condition in the definition of $\mathcal{S}_j$; that is, we replace $\mathcal{S}'_j$ with $\mathcal{S}''_j$, the set of $(y_1, \ldots, y_n)$ satisfying only the two conditions:

(i) $|y_1| \geq \cdots \geq |y_j|$ and $|y_{j+1}| \leq \cdots \leq |y_n|$,
(ii) $|\sum_{i=1}^{n} K^{(1/d)} y_i| \geq L$.

This yields the bound

$$p_d(n, L, j) \leq \frac{K^n}{j!(n-j)!} P\left(\left|\sum_{i=1}^{n} W_i\right| \geq \frac{L}{K^{(1/d)}}\right)$$

$$\leq \frac{K^n}{j!(n-j)!} P\left(\sum_{i=1}^{n} |W_i| \geq \frac{L}{K^{(1/d)}}\right). \qquad \square$$

PROOF OF PROPOSITION 2.3. From (3) and (4) we have

$$\tau_d(L) \leq \sum_{n=1}^{\infty} \sum_{j=0}^{n} \frac{K^n}{j!(n-j)!} P(|W_1 + \cdots + W_n| \geq c_2 L)$$

$$= \sum_{n=1}^{\infty} \frac{(2K)^n}{n!} P(|W_1 + \cdots + W_n| \geq c_2 L)$$

(10)
$$\leq e^{2K} \sum_{n=0}^{\infty} e^{-2K} \frac{(2K)^n}{n!} P(|W_1| + \cdots + |W_n| \geq c_2 L)$$

$$= e^{2K} P(U \geq c_2 L).$$

The last equality follows by taking $U$ to be a compound Poisson random variable $|W_1| + \cdots + |W_N|$ (where $N$ is independent of $|W_1|, |W_2|, \ldots$ and is Poisson with mean $2K$). By a standard calculation,

(11) $$E(e^{rU}) = e^{2KE(e^{r|W_1|} - 1)}.$$

This completes the proof of Proposition 2.3. □

PROOF OF THEOREM 1.1. Taking $n \approx \frac{L}{(\log L)^{1/d}}$ in (2) we get the lower bound in (1).

For the upper bound we use large deviation bounds on $P(U \geq c_2 L)$:

$$P(U \geq c_2 L) \leq \exp\left\{\inf_{r>0}\{\log(E(e^{rU})) - rc_2 L\}\right\}$$

$$(12) \qquad = \exp\left\{\inf_{r>0}\{2KE(e^{r|W_1|} - 1) - rc_2 L\}\right\}$$

$$\leq \exp\left\{\inf_{r>0}\{2K'(e^{cr^{d/(d-1)}} - rc_2 L)\}\right\}.$$

The last inequality follows from the fact that $E(e^{r|W_1|}) \leq e^{cr^{d/(d-1)}}$ for large $r$ (as can easily be shown). Taking $r \approx \alpha(\log L)^{(d-1)/d}$ for an appropriate constant $\alpha$ we get $\tau_d(L) \leq e^{-C_2 L(\log L)^{(d-1)/d}}$ for large $L$. □

PROOF OF THEOREM 1.2. The lower bound follows from Proposition 2.1 by taking $L=0$. For the upper bound, we note that in order for a path with the required properties to exist, there has to be a path from the origin, touching at least $\lfloor \frac{n}{2} \rfloor$ points with at most one change in direction. Using Proposition 2.2, this leads to

$$\rho_d(n) \leq \sum_{j=0}^{\lfloor n/2 \rfloor} p_d\left(\left\lfloor \frac{n}{2} \right\rfloor, 0, j\right)$$

$$\leq \sum_{j=0}^{\lfloor n/2 \rfloor} \frac{K_d^{\lfloor n/2 \rfloor}}{j!(\lfloor n/2 \rfloor - 1)!}$$

$$= \frac{(2K_d)^{\lfloor n/2 \rfloor}}{\lfloor n/2 \rfloor!}$$

$$\approx e^{(n/2)\log(2K_d) - (n/2)(\log n - \log 2 - 1)}$$

$$\leq e^{-c_2 n \log n},$$

for large $n$, proving the result. □

**3. Proof of Theorem 1.3.** We shall need the following simple fact. Define

$$L_d(a, b, y) = V(B(s, a) \cap B(t, b)),$$

where $s$ and $t$ are two points in $\mathbb{R}^d$ at distance $y$.

PROPOSITION 3.1. *For each fixed $y \geq 1$ we have*
$$\frac{L_d(1,y,y)}{\pi_d} \to 0,$$
as $d \to \infty$.

PROOF. Note that $L_d(1,y,y)$ is a lens contained in a cylinder of radius $r = r(y) < 1$ and height 1. The volume of this cylinder is $\pi_{d-1} r^d = o(\pi_d)$. □

PROOF OF THEOREM 1.3. For notational convenience, we write $B_{i,j}$ for $B(s_i, |x_j|)$, where as before $x_j = s_j - s_{j-1}$, for $j = 1, 2, \ldots$.

The set of points $(s_1, \ldots, s_k) \in (\mathbb{R}^d)^k$ satisfying $|x_1| \geq |x_2| \geq \cdots \geq |x_k|$ is denoted by $\mathcal{U}_k$. Furthermore, the subset of $\mathcal{U}_k$ which in addition satisfies
$$s_i \notin \bigcup_{m=1}^{i} B_{m-1,m}$$
for $i = 2, \ldots, j$ is denoted by $\mathcal{U}'_j$, for $j = 2, \ldots, k$. For convenience we define $\mathcal{U}'_1 = \mathcal{U}_k$. As in Section 2 we may now write

$$g_d(k) = \int_{s=(s_1,\ldots,s_k)\in\mathcal{U}'_k} e^{-V(B_{0,1})} e^{-V(B_{1,2}\setminus B_{0,1})} e^{-V(B_{2,3}\setminus(B_{0,1}\cup B_{1,2}))}$$
$$\times \cdots \times e^{-V(B_{k-1,k}\cup B_{k,k}\setminus \bigcup_{i=1}^{k-1} B_{i-1,i})} ds.$$

We now first estimate this from below, and after that we show that the error we make by doing this, tends to zero when the dimension tends to infinity. The first step is to replace the volumes in the exponents by the volume we would get without subtracting anything from the first set mentioned in each exponent. Thus

$$(13) \quad g_d(k) \geq \int_{s=(s_1,\ldots,s_k)\in\mathcal{U}'_k} e^{-\pi_d |x_1|^d} e^{-\pi_d |x_2|^d} \cdots e^{-\pi_d |x_{k-1}|^d} e^{-2\pi_d |x_k|^d} ds.$$

Writing the integrand in this formula as $W(s)$, the next step is to rewrite this as

$$(14) \quad \int_{\mathcal{U}_k} W(s) \, ds - \sum_{j=2}^{k} \int_{\mathcal{U}'_{j-1}\setminus\mathcal{U}'_j} W(s) \, ds,$$

that is, we have an integral over $\mathcal{U}_k$ as our leading term, and subtract from this the integral over those sequences $s$ which have a first index $j$ for which $s_j$ falls into a previous ball, $j = 2, \ldots, k$.

To compute the leading term, note that a simple change of variables gives that

$$\int_{\mathcal{U}_k} W(s)\,ds = \int_{0 \le y_k \le y_{k-1} \le \cdots \le y_1} e^{-y_1} \cdots e^{-y_{k-1}} e^{-2y_k}\,dy_k \cdots dy_1. \tag{15}$$

This integral can be computed explicitly, but its value is most easily found and understood via a simple probabilistic interpretation. Indeed, the integral is equal to $1/2$ times the probability that $Y_1 \ge Y_2 \ge \cdots \ge Y_k$, where the $Y_i$'s are independent, where $Y_1, Y_2, \ldots, Y_{k-1}$ are standard exponentially distributed and where $Y_k$ has an exponential distribution with parameter 2. The probability that $Y_1, Y_2, \ldots, Y_{k-1}$ are ordered this way is just $1/(k-1)!$. Since $Y_k$ has the same distribution as the minimum of two independent exponentially distributed random variables, the probability that one of these two will be the smallest among the $k+1$ random variables in question is simply $2/(k+1)$. It follows that the integral is equal to

$$\frac{1}{2}\frac{1}{(k-1)!}\frac{2}{(k+1)} = \frac{k}{(k+1)!}.$$

Next we will show that the remaining terms in (14) tend to 0 as the dimension $d$ tends to infinity. For this we need to bound

$$\int_{\mathcal{U}'_{j-1} \setminus \mathcal{U}'_j} W(s)\,ds,$$

for $j = 2, \ldots, k$. For given $j$, this integral is over those sequences for which $j$ is the first index so that $s_j$ falls into a previous ball. Hence

$$\int_{\mathcal{U}'_{j-1} \setminus \mathcal{U}'_j} W(s)\,ds$$

$$\le \int_{s_1 \in \mathbb{R}^d} e^{-\pi_d |x_1|^d} \int_{s_2 \in B_{1,1}} e^{-\pi_d |x_2|^d} \cdots \int_{s_{j-1} \in B_{j-2,j-2}} e^{-\pi_d |x_{j-1}|^d}$$

$$\times \int_{s_j \in B_{j-1,j-1} \cap \bigcup_{i=1}^{j-1} B_{i-1,i}} e^{-\pi_d |x_j|^d}\,ds_j \cdots ds_1$$

$$\le \int_{s_1 \in \mathbb{R}^d} e^{-\pi_d |x_1|^d} \int_{s_2 \in B_{1,1}} e^{-\pi_d |x_2|^d} \cdots \int_{s_{j-1} \in B_{j-2,j-2}} e^{-\pi_d |x_{j-1}|^d}$$

$$\times \frac{\sum_{i=1}^{j-1} V(B_{j-1,j-1} \cap B_{i-1,i})}{V(B_{j-1,j-1})} \int_{s_j \in B_{j-1,j-1}} e^{-\pi_d |x_j|^d}\,ds_j \cdots ds_1,$$

where in the last inequality we have used the fact that if $f(|x|)$ is decreasing in $|x|$, then for $y$ with $|y| \ge 1$ it is the case that

$$\int_{x \in B(0,1) \cap B(y,|y|)} f(|x|)\,dx \le \frac{V(B(0,1) \cap B(y,|y|))}{V(B(0,1))} \int_{x \in B(0,1)} f(|x|)\,dx.$$

The volume in the numerator can be estimated by the volume of the largest possible intersection, which is the one with the first ball in the situation that $s_{j-1}$ lies on its boundary:

$$\frac{\sum_{i=1}^{j-1} V(B_{j-1,j-1} \cap B_{i-1,i})}{V(B_{j-1,j-1})} \leq (j-1)\frac{L_d(|x_{j-1}|, |x_1|, |x_1|)}{V(B(0, |x_{j-1}|))}$$

$$= \frac{(j-1)}{\pi_d} L_d\left(1, \left|\frac{x_1}{x_{j-1}}\right|, \left|\frac{x_1}{x_{j-1}}\right|\right).$$

We now continue our estimate as follows, using the change of variables $y_i = \pi_d |x_i|^d$:

$$\int_{\mathcal{U}'_{j-1} \setminus \mathcal{U}'_j} W(s)\, ds$$

$$\leq \int_{|x_1|>0} e^{-\pi_d|x_1|^d} d\pi_d|x_1|^{d-1} \int_{|x_2|\leq|x_1|} e^{-\pi_d|x_2|^d} d\pi_d|x_2|^{d-1} \cdots$$

$$\times \int_{|x_{j-1}|\leq|x_{j-2}|} e^{-\pi_d|x_{j-1}|^d} d\pi_d|x_{j-1}|^{d-1} \frac{(j-1)}{\pi_d} L_d\left(1, \left|\frac{x_1}{x_{j-1}}\right|, \left|\frac{x_1}{x_{j-1}}\right|\right)$$

$$\times \int_{|x_j|\leq|x_{j-1}|} e^{-\pi_d|x_j|^d} d\pi_d|x_j|^{d-1} d|x_j|\cdots d|x_2|\, d|x_1|$$

$$= \int_{0\leq y_j\leq\cdots\leq y_1} \frac{(j-1)}{\pi_d} L_d\left(1, \left(\frac{y_1}{y_{j-1}}\right)^{1/d}, \left(\frac{y_1}{y_{j-1}}\right)^{1/d}\right)$$

$$\times e^{-\sum_{i=1}^{j} y_i} dy_j \cdots dy_1$$

$$\leq \int_{0\leq y_j\leq\cdots\leq y_1} \frac{(j-1)}{\pi_d} L_d\left(1, \frac{y_1}{y_{j-1}}, \frac{y_1}{y_{j-1}}\right) e^{-\sum_{i=1}^{j} y_i} dy_j \cdots dy_1,$$

which tends to zero according to Proposition 3.1 and dominated convergence.

It follows from all this that

$$\lim_{d\to\infty} g_d(k) \geq \frac{k}{(k+1)!}.$$

The final observation is the following: since the sum over $k$ of $k/(k+1)!$ is equal to 1, the inequality is in fact an equality, and we conclude that

$$\lim_{d\to\infty} g_d(k) = \frac{k}{(k+1)!}, \qquad k = 1, 2, \ldots,$$

proving the theorem. $\square$

**Acknowledgments.** We thank Charles Newman for useful conversations. We also thank an anonymous referee for his or her very careful reading of the manuscript, which greatly improved this paper.

I. KOZAKOVA  
R. MEESTER  
DEPARTMENT OF MATHEMATICS  
VRIJE UNIVERSITEIT AMSTERDAM  
DE BOELELAAN 1081  
1081 HV AMSTERDAM  
THE NETHERLANDS  
E-MAIL: rmeester@cs.vu.nl

S. NANDA  
TIFR CENTRE  
P.O. BOX 1234  
IISC CAMPUS  
BANGALORE 560012  
INDIA  
E-MAIL: nanda@math.tifrbng.res.in